\documentclass[12pt]{amsart}
\begin{document}
\title{Symmetry via Lie algebra cohomology}
\author{Michael Eastwood}
\address{Mathematical Sciences Institute\\
Australian National University\\ ACT 0200\\ Australia}
\email{meastwoo@member.ams.org}
\thanks{The author is supported by the Australian Research Council.}
\begin{abstract} The Killing operator on a Riemannian manifold is a linear 
differential operator on vector fields whose kernel provides the infinitesimal
Riemannian symmetries. The Killing operator is best understood in terms of its 
prolongation, which entails some simple tensor identities. These simple 
identities can be viewed as arising from the identification of certain Lie
algebra cohomologies. The point is that this case provides a model for 
more complicated operators similarly concerned with symmetry.
\end{abstract}
\maketitle

\vspace{-20pt}
\section{Disclaimer}
The results in this article are not widely known but are implicitly already
contained in~\cite{bceg,cd,css4}, for example. The object of this short
exposition is to introduce the method, by means of familiar examples, to a
wider audience.

\section{Notation}
The notation in this article follows the standard index conventions of
differential geometry. Precisely, we shall follow Penrose's abstract index
notation~\cite{OT} in which tensors are systematically adorned with indices to
specify their type. For example, vector fields are denoted with an upper index
$X^a$ whilst $2$-forms have $2$ lower indices~$\omega_{ab}$. The natural
contraction between them is denoted by repeating an index $X^a\omega_{ab}$ in
accordance with the Einstein summation convention. Round brackets are used to
denote symmetrisation over the indices they enclose whilst square brackets are
used to denote skewing, e.g.
$$\textstyle\psi_{[abc]d}=\frac16[\psi_{abcd}+\psi_{bcad}+\psi_{cabd}
-\psi_{bacd}-\psi_{acbd}-\psi_{cbad}].$$

\section{The Levi-Civita connection}
Suppose $g_{ab}$ is a Riemannian metric. The Levi-Civita connection $\nabla_a$
associated with $g_{ab}$ is characterised by the following well-known
properties
\begin{itemize}
\item $\nabla_a$ is torsion-free,
\item $\nabla_ag_{bc}=0$.
\end{itemize}
Its existence and uniqueness boils down to a tensor identity as follows. 
Choose~$D_a$, any torsion-free connection. Any other must be of the form 
$$\nabla_a\phi_b=D_a\phi_b-\Gamma_{ab}{}^c\phi_c$$
for some tensor $\Gamma_{ab}{}^c=\Gamma_{(ab)}{}^c$ and then $\nabla_ag_{bc}=0$
if and only if 
$$0=D_ag_{bc}-\Gamma_{ab}{}^dg_{dc}-\Gamma_{ac}{}^dg_{bd}=
D_ag_{bc}-\Gamma_{abc}-\Gamma_{acb},$$
where we are using the metric $g_{ab}$ to `lower indices' in the usual fashion.
These are two conditions on $\Gamma_{abc}$, namely
$$\textstyle\Gamma_{[ab]c}=0\qquad\mbox{and}\qquad
\Gamma_{a(bc)}=\frac12D_ag_{bc}$$
that always have a unique solution. To see this, note that the general 
solution of the second equation has the form
$$\textstyle\Gamma_{abc}=\frac12D_ag_{bc}-K_{abc},\quad
\mbox{where }K_{abc}=K_{a[bc]}.$$
Having done this, the first equation reads
$$\textstyle K_{[ab]c}=\frac12D_{[a}g_{b]c},$$
which always has a unique solution owing to the tensor isomorphism 
\begin{equation}\label{key}\framebox{$\begin{array}{ccl}
\Lambda^1\otimes\Lambda^2&
\begin{array}[b]l\simeq\\[-9pt] \longrightarrow\end{array}&
\Lambda^2\otimes\Lambda^1\\[3pt]
K_{abc}=K_{a[bc]}&\longmapsto&K_{[ab]c}\,,
\end{array}$}\end{equation}
where $\Lambda^p$ denotes the bundle of $p$-forms. This isomorphism is typical
of the tensor identities to be explained in this article by means of Lie
algebra cohomology.

\section{The Killing operator}
A vector field $X^a$ on a Riemannian manifold with metric $g_{ab}$ is said to
be a {\em Killing field\/} if and only if ${\mathcal{L}}_Xg_{ab}=0$, where
${\mathcal{L}}_X$ is the Lie derivative along~$X^a$. The geometric
interpretation of Lie derivative means that the flow of $X^a$ is an isometry. 
Thus, a Killing field is an {\em infinitesimal symmetry\/} in the context of 
Riemannian geometry.

It is useful to regard the Killing equation ${\mathcal{L}}_Xg_{ab}=0$ as a
linear partial differential equation on the vector field $X^a$ as follows. For
any torsion-free connection~$\nabla_a$,
$${\mathcal{L}}_X\phi_{b}=X^a\nabla_a\phi_b+\phi_a\nabla_bX^a$$
so, if we use the Levi-Civita connection for $g_{ab}$, then
$$\begin{array}{rcl}{\mathcal{L}}_Xg_{bc}&=&
X^a\nabla_ag_{bc}+g_{ac}\nabla_bX^a+g_{ba}\nabla_cX^a\\[3pt]
&=&\nabla_bX_c+\nabla_cX_b.
\end{array}$$
Hence, the Killing fields $X^a$ make up the kernel of the 
{\em Killing operator\/}:-- 
$$\begin{array}{rcccl}\mbox{Tangent bundle}&
\begin{array}[b]l\simeq\\[-9pt] \longrightarrow\end{array}&
\Lambda^1&\longrightarrow&\textstyle\bigodot^2\!\Lambda^1\\[3pt]
X^a&\longmapsto&X_a&\longmapsto&\nabla_{(a}X_{b)}\,.
\end{array}$$

\section{Prolongation of the Killing operator}\label{prolongedkilling}
For any torsion-free connection $\nabla_a$, the equation $\nabla_{(a}X_{b)}=0$
may be understood as follows. Certainly, we may rewrite it as
\begin{equation}\label{killing}
\nabla_aX_b=K_{ab},\quad\mbox{where $K_{ab}$ is skew}.\end{equation}
In this case $\nabla_{[a}K_{bc]}=0$, a condition which we may rewrite as 
$$\nabla_aK_{bc}=\nabla_cK_{ba}-\nabla_bK_{ca}$$
and substitute from (\ref{killing}) to conclude, as a differential 
consequence, that
$$\nabla_aK_{bc}=\nabla_c\nabla_bX_a-\nabla_b\nabla_cX_a=R_{bc}{}^d{}_aX_d,$$
where $R_{ab}{}^c{}_d$ is the curvature of $\nabla_a$ characterised by 
$$[\nabla_a\nabla_b-\nabla_b\nabla_a]X^c=R_{ab}{}^c{}_dX^d.$$
Therefore,
$$\nabla_{(a}X_{b)}=0\iff\framebox{$\begin{array}{rcl}
\nabla_aX_b&=&K_{ab}\\
\nabla_aK_{bc}&=&R_{bc}{}^d{}_aX_d
\end{array}$}$$
In other words, Killing fields are in $1$--$1$ correspondence with covariant
constant sections of the vector bundle
${\mathbb{T}}\equiv\Lambda^1\oplus\Lambda^2$ equipped with the connection
\begin{equation}\label{tractorconnection}
{\mathbb{T}}\ni\left[\begin{array}cX_b\\ K_{bc}\end{array}\right]
\stackrel{\mbox{\small$\nabla_a$}}{\longmapsto}
\left[\begin{array}c\nabla_aX_b-K_{ab}\\ 
\nabla_aK_{bc}-R_{bc}{}^d{}_aX_d\end{array}\right]\in
\Lambda^1\otimes{\mathbb{T}}.\end{equation}
At this point, we may use the standard theory of vector bundles with connection
to investigate Killing fields. In particular, it is immediately clear that the
Killing fields on a connected manifold form a vector space whose dimension is 
bounded by the rank of ${\mathbb{T}}$, namely $n(n+1)/2$. 

\section{The Killing operator in flat space}
Be that as it may, suppose ask only about the Killing operator on 
flat space. It is easily verified in this case that the connection 
(\ref{tractorconnection}) is flat (and, in fact, the same is true on any 
constant curvature space). Therefore, we may couple the de~Rham sequence with 
(\ref{tractorconnection}) to obtain a locally exact complex
$${\mathbb{T}}\stackrel{\nabla}{\longrightarrow}
\Lambda^1\otimes{\mathbb{T}}\stackrel{\nabla}{\longrightarrow}
\Lambda^2\otimes{\mathbb{T}}\stackrel{\nabla}{\longrightarrow}
\Lambda^3\otimes{\mathbb{T}}\stackrel{\nabla}{\longrightarrow}\cdots$$
and, at this point, the isomorphism (\ref{key}) re-emerges! Specifically, in 
the absence of the curvature term (\ref{tractorconnection}) may be written as
$$\left[\begin{array}cX_b\\ K_{bc}\end{array}\right]
\stackrel{\mbox{\small$\nabla_a$}}{\longmapsto}
\left[\begin{array}c\nabla_aX_b\\ 
\nabla_aK_{bc}\end{array}\right]
-\partial\!\left[\begin{array}cX_b\\ K_{bc}
\end{array}\right],\enskip\mbox{where }
\partial\!\left[\begin{array}cX_b\\ 
K_{bc}\end{array}\right]
=\left[\begin{array}cK_{ab}\\ 0
\end{array}\right].$$
The homomorphism $\partial:{\mathbb{T}}\to\Lambda^1\otimes{\mathbb{T}}$ induces
$\partial:\Lambda^p\otimes{\mathbb{T}}\to\Lambda^{p+1}\otimes{\mathbb{T}}$ by
$\partial(\omega\otimes X)=\omega\wedge\partial X$ and we obtain a complex
\begin{equation}\label{coupledtractors}
\makebox[200pt]{$\begin{array}{ccccccccccc}
0&\to&{\mathbb{T}}
&\stackrel{\partial}{\longrightarrow}&\Lambda^1\otimes{\mathbb{T}}
&\stackrel{\partial}{\longrightarrow}&\Lambda^2\otimes{\mathbb{T}}
&\stackrel{\partial}{\longrightarrow}&\Lambda^3\otimes{\mathbb{T}}
&\stackrel{\partial}{\longrightarrow}&\cdots\\
&&\|&&\|&&\|&&\|\\
&&\Lambda^1&&\Lambda^1\otimes\Lambda^1&&
\Lambda^2\otimes\Lambda^1
&&\Lambda^3\otimes\Lambda^1&\!{}\!&\cdots\\
&&\oplus&\begin{picture}(0,0)
\put(-10,-10){\vector(1,1){22}}
\end{picture}&\oplus&\begin{picture}(0,0)
\put(-10,-10){\vector(1,1){22}}
\put(-20,5){\framebox{\tiny NB}}
\end{picture}&\oplus&\begin{picture}(0,0)
\put(-10,-10){\vector(1,1){22}}
\end{picture}&\oplus&\begin{picture}(0,0)
\put(-10,-10){\vector(1,1){22}}
\end{picture}\\
&&\Lambda^2&&\Lambda^1\otimes\Lambda^2&&\Lambda^2\otimes\Lambda^2
&&\Lambda^3\otimes\Lambda^2&\!{}\!&\cdots
\end{array}$}\end{equation}
in which 
$\partial:\Lambda^1\otimes{\mathbb{T}}\to\Lambda^2\otimes{\mathbb{T}}$ is 
carried by the isomorphism~(\ref{key}). More generally, we can ask about the 
cohomology of the complex $(\Lambda^\bullet\otimes{\mathbb{T}},\partial)$ and 
conclude, by inspection, that
$$\begin{array}{rcl}
H^0(\Lambda^\bullet\otimes{\mathbb{T}},\partial)&=&\{X_a\}\\
H^1(\Lambda^\bullet\otimes{\mathbb{T}},\partial)&=&\{X_{ab}=X_{(ab)}\}\\
H^2(\Lambda^\bullet\otimes{\mathbb{T}},\partial)&=&
\{K_{abcd}=K_{[ab][cd]}\mbox{ s.t.\ }K_{[abc]d}=0\}\\
H^3(\Lambda^\bullet\otimes{\mathbb{T}},\partial)&=&
\{K_{abcde}=K_{[abc][de]}\mbox{ s.t.\ }K_{[abcd]e}=0\}\\
H^4(\Lambda^\bullet\otimes{\mathbb{T}},\partial)&=&
\{K_{abcdef}=K_{[abcd][ef]}\mbox{ s.t.\ }K_{[abcde]f}=0\}\\
\vdots&\vdots&\vdots
\end{array}$$
recognising that each of these bundles is an irreducible tensor bundle, 
which we may write as Young diagrams~\cite{fh}

\vspace*{-25pt}
\begin{equation}\label{kostant}H^0=
\begin{picture}(6,6)
\put(0,0){\line(1,0){6}}
\put(0,6){\line(1,0){6}}
\put(0,0){\line(0,1){6}}
\put(6,0){\line(0,1){6}}
\end{picture}\quad
H^1=
\begin{picture}(12,6)
\put(0,0){\line(1,0){12}}
\put(0,6){\line(1,0){12}}
\put(0,0){\line(0,1){6}}
\put(6,0){\line(0,1){6}}
\put(12,0){\line(0,1){6}}
\end{picture}\quad
H^2=
\begin{picture}(12,12)(0,6)
\put(0,0){\line(1,0){12}}
\put(0,6){\line(1,0){12}}
\put(0,12){\line(1,0){12}}
\put(0,0){\line(0,1){12}}
\put(6,0){\line(0,1){12}}
\put(12,0){\line(0,1){12}}
\end{picture}\quad
H^3=
\begin{picture}(12,18)(0,12)
\put(0,0){\line(1,0){6}}
\put(0,6){\line(1,0){12}}
\put(0,12){\line(1,0){12}}
\put(0,18){\line(1,0){12}}
\put(0,0){\line(0,1){18}}
\put(6,0){\line(0,1){18}}
\put(12,6){\line(0,1){12}}
\end{picture}\quad
H^4=
\begin{picture}(12,24)(0,18)
\put(0,0){\line(1,0){6}}
\put(0,6){\line(1,0){6}}
\put(0,12){\line(1,0){12}}
\put(0,18){\line(1,0){12}}
\put(0,24){\line(1,0){12}}
\put(0,0){\line(0,1){24}}
\put(6,0){\line(0,1){24}}
\put(12,12){\line(0,1){12}}
\end{picture}
\quad\cdots\end{equation}
\vspace*{3pt}

\noindent Readers may notice that
$H^2(\Lambda^\bullet\otimes{\mathbb{T}},\partial)$ is the natural location for
the Riemann curvature tensor and that 
$H^3(\Lambda^\bullet\otimes{\mathbb{T}},\partial)$ is the natural location for 
the Bianchi identity. These observations are more fully explained 
in~\cite{imanotes}. Here, suffice it to observe that a simple diagram chase on 
(\ref{coupledtractors}) reveals a locally exact complex

\vspace*{-25pt}
\begin{equation}\label{simpleBGG}\begin{picture}(6,6)
\put(0,0){\line(1,0){6}}
\put(0,6){\line(1,0){6}}
\put(0,0){\line(0,1){6}}
\put(6,0){\line(0,1){6}}
\end{picture}\xrightarrow{\,\nabla\,}
\begin{picture}(12,6)
\put(0,0){\line(1,0){12}}
\put(0,6){\line(1,0){12}}
\put(0,0){\line(0,1){6}}
\put(6,0){\line(0,1){6}}
\put(12,0){\line(0,1){6}}
\end{picture}\xrightarrow{\,\nabla^{(2)}\,}
\begin{picture}(12,12)(0,6)
\put(0,0){\line(1,0){12}}
\put(0,6){\line(1,0){12}}
\put(0,12){\line(1,0){12}}
\put(0,0){\line(0,1){12}}
\put(6,0){\line(0,1){12}}
\put(12,0){\line(0,1){12}}
\end{picture}\xrightarrow{\,\nabla\,}
\begin{picture}(12,18)(0,12)
\put(0,0){\line(1,0){6}}
\put(0,6){\line(1,0){12}}
\put(0,12){\line(1,0){12}}
\put(0,18){\line(1,0){12}}
\put(0,0){\line(0,1){18}}
\put(6,0){\line(0,1){18}}
\put(12,6){\line(0,1){12}}
\end{picture}\xrightarrow{\,\nabla\,}
\begin{picture}(12,24)(0,18)
\put(0,0){\line(1,0){6}}
\put(0,6){\line(1,0){6}}
\put(0,12){\line(1,0){12}}
\put(0,18){\line(1,0){12}}
\put(0,24){\line(1,0){12}}
\put(0,0){\line(0,1){24}}
\put(6,0){\line(0,1){24}}
\put(12,12){\line(0,1){12}}
\end{picture}\xrightarrow{\,\nabla\,}\cdots\end{equation}
\vspace*{1pt}

\noindent and, in particular, an identification of the range of the Killing 
operator in flat space as follows.

\smallskip\noindent{\bf Theorem}\quad {\em Suppose 
$U$ is an open subset of ${\mathbb{R}}^n$ with $H^1(U,{\mathbb{R}})=0$. Then a 
symmetric tensor $\omega_{ab}$ on $U$ is of the form $\nabla_{(a}X_{b)}$ for 
some $X_a$ on $U$ if and only if}
$$\nabla_a\nabla_c\omega_{bd}-\nabla_b\nabla_c\omega_{ad}
-\nabla_a\nabla_d\omega_{bc}+\nabla_b\nabla_d\omega_{ac}=0.$$

\section{Higher Killing operators}
So far, we have not seen any Lie algebra cohomology, although it is lurking in
the background. The identifications (\ref{kostant}) can be obtained by
elementary means. As soon as we consider more complicated operators, however,
then the corresponding identifications are not so obvious. A {\em Killing
tensor\/} of valence $\ell$ is a symmetric tensor field $X_{bc\cdots de}$ with
$\ell$ indices annihilated by the higher Killing operator
$$X_{bc\cdots de}\mapsto\nabla_{(a}X_{bc\cdots de)}.$$
Killing tensors induce conserved quantities along geodesics and arise naturally
in the theory of separation of variables. The higher Killing operators may be
prolonged along the lines explained in~\S\ref{prolongedkilling}. The details
are more complicated and this is where Lie algebra cohomology comes to the
fore. Without going into details, the prolonged bundle
$${\mathbb{T}}=\Lambda^1\oplus\Lambda^2=\begin{picture}(6,6)
\put(0,0){\line(1,0){6}}
\put(0,6){\line(1,0){6}}
\put(0,0){\line(0,1){6}}
\put(6,0){\line(0,1){6}}
\end{picture}\oplus
\begin{picture}(6,12)(0,6)
\put(0,0){\line(1,0){6}}
\put(0,6){\line(1,0){6}}
\put(0,12){\line(1,0){6}}
\put(0,0){\line(0,1){12}}
\put(6,0){\line(0,1){12}}
\end{picture}$$ 
that we saw in \S\ref{prolongedkilling} should be replaced by
$${\mathbb{T}}^\ell={\mathbb{T}}_0^\ell\oplus\cdots\oplus{\mathbb{T}}_\ell^\ell
=\begin{picture}(42,6)
\put(0,0){\line(1,0){42}}
\put(0,6){\line(1,0){42}}
\put(0,0){\line(0,1){6}}
\put(6,0){\line(0,1){6}}
\put(12,0){\line(0,1){6}}
\put(18,0){\line(0,1){6}}
\put(36,0){\line(0,1){6}}
\put(42,0){\line(0,1){6}}
\put(23,3){\makebox(0,0){.}}
\put(27,3){\makebox(0,0){.}}
\put(31,3){\makebox(0,0){.}}
\put(0,-8){\line(0,1){4}}
\put(42,-8){\line(0,1){4}}
\put(21,-5){\makebox(0,0){\tiny$\ell$ boxes}}
\put(6,-6){\vector(-1,0){6}}
\put(36,-6){\vector(1,0){6}}
\end{picture}\oplus\begin{picture}(42,12)(0,6)
\put(0,0){\line(1,0){6}}
\put(0,6){\line(1,0){42}}
\put(0,12){\line(1,0){42}}
\put(0,0){\line(0,1){12}}
\put(6,0){\line(0,1){12}}
\put(12,6){\line(0,1){6}}
\put(18,6){\line(0,1){6}}
\put(36,6){\line(0,1){6}}
\put(42,6){\line(0,1){6}}
\put(23,9){\makebox(0,0){.}}
\put(27,9){\makebox(0,0){.}}
\put(31,9){\makebox(0,0){.}}
\end{picture}\oplus\begin{picture}(42,12)(0,6)
\put(0,0){\line(1,0){12}}
\put(0,6){\line(1,0){42}}
\put(0,12){\line(1,0){42}}
\put(0,0){\line(0,1){12}}
\put(6,0){\line(0,1){12}}
\put(12,0){\line(0,1){12}}
\put(18,6){\line(0,1){6}}
\put(36,6){\line(0,1){6}}
\put(42,6){\line(0,1){6}}
\put(23,9){\makebox(0,0){.}}
\put(27,9){\makebox(0,0){.}}
\put(31,9){\makebox(0,0){.}}
\end{picture}\oplus\cdots\oplus\begin{picture}(42,12)(0,6)
\put(0,0){\line(1,0){42}}
\put(0,6){\line(1,0){42}}
\put(0,12){\line(1,0){42}}
\put(0,0){\line(0,1){12}}
\put(6,0){\line(0,1){12}}
\put(12,0){\line(0,1){12}}
\put(18,0){\line(0,1){12}}
\put(36,0){\line(0,1){12}}
\put(42,0){\line(0,1){12}}
\put(23,3){\makebox(0,0){.}}
\put(27,3){\makebox(0,0){.}}
\put(31,3){\makebox(0,0){.}}
\put(23,9){\makebox(0,0){.}}
\put(27,9){\makebox(0,0){.}}
\put(31,9){\makebox(0,0){.}}
\end{picture}\,,$$
realised as
$$\left[\begin{array}l
X_{bc\cdots de}=X_{(bc\cdots de)}\\
K_{pbc\cdots de}=K_{p(bc\cdots de)}\mbox{ s.t.\ }K_{(pbc\cdots de)}=0\\
K_{pqbc\cdots de}^{\prime}=K_{(pq)(bc\cdots de)}^{\prime}
\mbox{ s.t.\ }K_{p(qbc\cdots de)}^{\prime}=0\\
K_{pqrbc\cdots de}^{\prime\prime}=K_{(pqr)(bc\cdots de)}^{\prime\prime}
\mbox{ s.t.\ }K_{pq(rbc\cdots de)}^{\prime\prime}=0\\
\hspace{40pt}\vdots\\
K_{pq\cdots rsbc\cdots de}^{\prime\prime\cdots\prime}=
K_{(pq\cdots rs)(bc\cdots de)}^{\prime\prime\cdots\prime}\mbox{ s.t.\ }
K_{pq\cdots r(sbc\cdots de)}^{\prime\prime\cdots\prime}=0
\end{array}\right]$$
with $\partial:{\mathbb{T}}^\ell\to\Lambda^1\otimes{\mathbb{T}}^\ell$
defined by
$$\partial\left[\begin{array}c
X_{bc\cdots de}\\
K_{pbc\cdots de}\\
K_{pqbc\cdots de}^{\prime}\\
K_{pqrbc\cdots de}^{\prime\prime}\\
\vdots\\
K_{pq\cdots rsbc\cdots de}^{\prime\prime\cdots\prime}
\end{array}\right]=
\left[\begin{array}c
K_{abc\cdots de}\\
K_{apbc\cdots de}^\prime\\
K_{apqbc\cdots de}^{\prime\prime}\\
\vdots\\
K_{ap\cdots qrbc\cdots de}^{\prime\prime\cdots\prime}\\
0
\end{array}\right].$$
The identifications generalising (\ref{kostant}) are as follows.
\begin{equation}\label{gen1}
H^0(\Lambda^\bullet\otimes{\mathbb{T}}^\ell,\partial)=
\begin{picture}(42,6)
\put(0,0){\line(1,0){42}}
\put(0,6){\line(1,0){42}}
\put(0,0){\line(0,1){6}}
\put(6,0){\line(0,1){6}}
\put(12,0){\line(0,1){6}}
\put(18,0){\line(0,1){6}}
\put(36,0){\line(0,1){6}}
\put(42,0){\line(0,1){6}}
\put(23,3){\makebox(0,0){.}}
\put(27,3){\makebox(0,0){.}}
\put(31,3){\makebox(0,0){.}}
\put(0,-8){\line(0,1){4}}
\put(42,-8){\line(0,1){4}}
\put(21,-5){\makebox(0,0){\tiny$\ell$ boxes}}
\put(6,-6){\vector(-1,0){6}}
\put(36,-6){\vector(1,0){6}}
\end{picture}\qquad
H^1(\Lambda^\bullet\otimes{\mathbb{T}}^\ell,\partial)=
\begin{picture}(48,6)
\put(0,0){\line(1,0){48}}
\put(0,6){\line(1,0){48}}
\put(0,0){\line(0,1){6}}
\put(6,0){\line(0,1){6}}
\put(12,0){\line(0,1){6}}
\put(18,0){\line(0,1){6}}
\put(24,0){\line(0,1){6}}
\put(42,0){\line(0,1){6}}
\put(48,0){\line(0,1){6}}
\put(29,3){\makebox(0,0){.}}
\put(33,3){\makebox(0,0){.}}
\put(37,3){\makebox(0,0){.}}
\put(0,-8){\line(0,1){4}}
\put(48,-8){\line(0,1){4}}
\put(24,-5){\makebox(0,0){\tiny$\ell\!\!+\!\!1$ boxes}}
\put(6,-6){\vector(-1,0){6}}
\put(42,-6){\vector(1,0){6}}
\end{picture}\end{equation}
and

\vspace*{-35pt}
\begin{equation}\label{gen2}\raisebox{20pt}{
\makebox[0pt]{$H^2(\Lambda^\bullet\otimes{\mathbb{T}}^\ell,\partial)=
\begin{picture}(48,12)(0,6)
\put(0,0){\line(1,0){48}}
\put(0,6){\line(1,0){48}}
\put(0,12){\line(1,0){48}}
\put(0,0){\line(0,1){12}}
\put(6,0){\line(0,1){12}}
\put(12,0){\line(0,1){12}}
\put(18,0){\line(0,1){12}}
\put(24,0){\line(0,1){12}}
\put(42,0){\line(0,1){12}}
\put(48,0){\line(0,1){12}}
\put(29,3){\makebox(0,0){.}}
\put(33,3){\makebox(0,0){.}}
\put(37,3){\makebox(0,0){.}}
\put(29,9){\makebox(0,0){.}}
\put(33,9){\makebox(0,0){.}}
\put(37,9){\makebox(0,0){.}}
\put(0,-8){\line(0,1){4}}
\put(48,-8){\line(0,1){4}}
\put(24,-5){\makebox(0,0){\tiny$\ell\!\!+\!\!1$ boxes}}
\put(6,-6){\vector(-1,0){6}}
\put(42,-6){\vector(1,0){6}}
\end{picture}\quad\enskip
H^p(\Lambda^\bullet\otimes{\mathbb{T}}^\ell,\partial)=
\begin{picture}(58,36)(0,30)
\put(0,0){\line(1,0){6}}
\put(0,6){\line(1,0){6}}
\put(0,24){\line(1,0){48}}
\put(0,30){\line(1,0){48}}
\put(0,36){\line(1,0){48}}
\put(0,0){\line(0,1){36}}
\put(6,0){\line(0,1){36}}
\put(12,24){\line(0,1){12}}
\put(18,24){\line(0,1){12}}
\put(24,24){\line(0,1){12}}
\put(42,24){\line(0,1){12}}
\put(48,24){\line(0,1){12}}
\put(29,27){\makebox(0,0){.}}
\put(33,27){\makebox(0,0){.}}
\put(37,27){\makebox(0,0){.}}
\put(29,33){\makebox(0,0){.}}
\put(33,33){\makebox(0,0){.}}
\put(37,33){\makebox(0,0){.}}
\put(3,11){\makebox(0,0){.}}
\put(3,15){\makebox(0,0){.}}
\put(3,19){\makebox(0,0){.}}
\put(52,0){\line(1,0){4}}
\put(52,36){\line(1,0){4}}
\put(55,18){\makebox(0,0){\small$p$}}
\put(54,12){\vector(0,-1){12}}
\put(54,24){\vector(0,1){12}}
\end{picture}
\mbox{ for }p\geq 2.$}}\end{equation}

The locally exact complex generalising (\ref{simpleBGG}) is

\vspace*{-15pt}
\begin{equation}\label{BGG}\begin{picture}(42,6)
\put(0,0){\line(1,0){42}}
\put(0,6){\line(1,0){42}}
\put(0,0){\line(0,1){6}}
\put(6,0){\line(0,1){6}}
\put(12,0){\line(0,1){6}}
\put(18,0){\line(0,1){6}}
\put(36,0){\line(0,1){6}}
\put(42,0){\line(0,1){6}}
\put(23,3){\makebox(0,0){.}}
\put(27,3){\makebox(0,0){.}}
\put(31,3){\makebox(0,0){.}}
\end{picture}\xrightarrow{\,\nabla\,}
\begin{picture}(48,6)
\put(0,0){\line(1,0){48}}
\put(0,6){\line(1,0){48}}
\put(0,0){\line(0,1){6}}
\put(6,0){\line(0,1){6}}
\put(12,0){\line(0,1){6}}
\put(18,0){\line(0,1){6}}
\put(24,0){\line(0,1){6}}
\put(42,0){\line(0,1){6}}
\put(48,0){\line(0,1){6}}
\put(29,3){\makebox(0,0){.}}
\put(33,3){\makebox(0,0){.}}
\put(37,3){\makebox(0,0){.}}
\end{picture}\xrightarrow{\,\nabla^{(\ell+1)}\,}
\begin{picture}(48,12)(0,6)
\put(0,0){\line(1,0){48}}
\put(0,6){\line(1,0){48}}
\put(0,12){\line(1,0){48}}
\put(0,0){\line(0,1){12}}
\put(6,0){\line(0,1){12}}
\put(12,0){\line(0,1){12}}
\put(18,0){\line(0,1){12}}
\put(24,0){\line(0,1){12}}
\put(42,0){\line(0,1){12}}
\put(48,0){\line(0,1){12}}
\put(29,3){\makebox(0,0){.}}
\put(33,3){\makebox(0,0){.}}
\put(37,3){\makebox(0,0){.}}
\put(29,9){\makebox(0,0){.}}
\put(33,9){\makebox(0,0){.}}
\put(37,9){\makebox(0,0){.}}
\end{picture}\xrightarrow{\,\nabla\,}
\begin{picture}(48,18)(0,12)
\put(0,0){\line(1,0){6}}
\put(0,6){\line(1,0){48}}
\put(0,12){\line(1,0){48}}
\put(0,18){\line(1,0){48}}
\put(0,0){\line(0,1){18}}
\put(6,0){\line(0,1){18}}
\put(12,6){\line(0,1){12}}
\put(18,6){\line(0,1){12}}
\put(24,6){\line(0,1){12}}
\put(42,6){\line(0,1){12}}
\put(48,6){\line(0,1){12}}
\put(29,9){\makebox(0,0){.}}
\put(33,9){\makebox(0,0){.}}
\put(37,9){\makebox(0,0){.}}
\put(29,15){\makebox(0,0){.}}
\put(33,15){\makebox(0,0){.}}
\put(37,15){\makebox(0,0){.}}
\end{picture}\xrightarrow{\,\nabla\,}\cdots\end{equation}

\vspace*{10pt}
\noindent where the first operator is the higher Killing operator. It is a 
special case of the Bernstein-Gelfand-Gelfand resolution~\cite{cd,css4}. 

\section{Tensor identities}
Be that as it may, the identifications of
$H^p(\Lambda^\bullet\otimes{\mathbb{T}}^\ell,\partial)$ claimed in the previous
section are not so easy and entail some tricky tensor identities. The natural
generalisation of~(\ref{key}), for example, follows by writing out the complex
$$0\to{\mathbb{T}}^\ell\stackrel{\partial}{\longrightarrow}
\Lambda^1\otimes{\mathbb{T}}^\ell\stackrel{\partial}{\longrightarrow}
\Lambda^2\otimes{\mathbb{T}}^\ell\stackrel{\partial}{\longrightarrow}
\Lambda^3\otimes{\mathbb{T}}^\ell\stackrel{\partial}{\longrightarrow}
\Lambda^4\otimes{\mathbb{T}}^\ell\stackrel{\partial}{\longrightarrow}
\cdots$$
as in (\ref{coupledtractors}) and pinning down the locations of the
cohomologies
$$\begin{array}{cccccccccc}
0&&\framebox{${\mathbb{T}}_0^\ell$}
&&\framebox{$\Lambda^1\otimes{\mathbb{T}}_0^\ell$}
&&\Lambda^2\otimes{\mathbb{T}}_0^\ell
&&\Lambda^3\otimes{\mathbb{T}}_0^\ell\\
\oplus&\nearrow&\oplus&\nearrow&\oplus&\nearrow&\oplus&\nearrow&\oplus
&\nearrow\\
0&&{\mathbb{T}}_1^\ell
&&\Lambda^1\otimes{\mathbb{T}}_1^\ell
&&\Lambda^2\otimes{\mathbb{T}}_1^\ell
&&\Lambda^3\otimes{\mathbb{T}}_1^\ell\\
\oplus&\nearrow&\oplus&\nearrow&\oplus&\nearrow&\oplus&\nearrow&\oplus
&\nearrow\\
\vdots&&\vdots&&\vdots&&\vdots&&\vdots\\
\oplus&\nearrow&\oplus&\nearrow&\oplus&\nearrow&\oplus&\nearrow&\oplus
&\nearrow\\
0&&{\mathbb{T}}_{\ell-2}^\ell
&&\Lambda^1\otimes{\mathbb{T}}_{\ell-2}^\ell\hspace*{-10pt}
&&\Lambda^2\otimes{\mathbb{T}}_{\ell-2}^\ell\hspace*{-10pt}
&&\Lambda^3\otimes{\mathbb{T}}_{\ell-2}^\ell\hspace*{-10pt}\\
\oplus&\nearrow&\oplus&\nearrow&\oplus&\nearrow&\oplus&\nearrow&\oplus
&\nearrow\\
0&&{\mathbb{T}}_{\ell-1}^\ell
&&\Lambda^1\otimes{\mathbb{T}}_{\ell-1}^\ell\hspace*{-10pt}
&&\Lambda^2\otimes{\mathbb{T}}_{\ell-1}^\ell\hspace*{-10pt}
&&\Lambda^3\otimes{\mathbb{T}}_{\ell-1}^\ell\hspace*{-10pt}\\
\oplus&\nearrow&\oplus&\nearrow&\oplus&\nearrow&\oplus&\nearrow&\oplus
&\nearrow\\
0&&{\mathbb{T}}_\ell^\ell
&&\Lambda^1\otimes{\mathbb{T}}_\ell^\ell
&&\framebox{$\Lambda^2\otimes{\mathbb{T}}_\ell^\ell$}
&&\framebox{$\Lambda^3\otimes{\mathbb{T}}_\ell^\ell$}
\end{array}$$
simply by the number of boxes involved to deduce that
$$\begin{array}r
0\to\Lambda^1\otimes\begin{picture}(42,12)(0,3)
\put(0,0){\line(1,0){42}}
\put(0,6){\line(1,0){42}}
\put(0,12){\line(1,0){42}}
\put(0,0){\line(0,1){12}}
\put(6,0){\line(0,1){12}}
\put(12,0){\line(0,1){12}}
\put(30,0){\line(0,1){12}}
\put(36,0){\line(0,1){12}}
\put(42,0){\line(0,1){12}}
\put(17,3){\makebox(0,0){.}}
\put(21,3){\makebox(0,0){.}}
\put(24,3){\makebox(0,0){.}}
\put(17,9){\makebox(0,0){.}}
\put(21,9){\makebox(0,0){.}}
\put(24,9){\makebox(0,0){.}}
\put(0,-8){\line(0,1){4}}
\put(42,-8){\line(0,1){4}}
\put(21,-5){\makebox(0,0){\scriptsize $\ell$}}
\put(14,-6){\vector(-1,0){14}}
\put(28,-6){\vector(1,0){14}}
\end{picture}\xrightarrow{\,\partial\,}
\Lambda^2\otimes\begin{picture}(42,12)(0,3)
\put(0,0){\line(1,0){36}}
\put(0,6){\line(1,0){42}}
\put(0,12){\line(1,0){42}}
\put(0,0){\line(0,1){12}}
\put(6,0){\line(0,1){12}}
\put(12,0){\line(0,1){12}}
\put(30,0){\line(0,1){12}}
\put(36,0){\line(0,1){12}}
\put(42,6){\line(0,1){6}}
\put(17,3){\makebox(0,0){.}}
\put(21,3){\makebox(0,0){.}}
\put(24,3){\makebox(0,0){.}}
\put(17,9){\makebox(0,0){.}}
\put(21,9){\makebox(0,0){.}}
\put(24,9){\makebox(0,0){.}}
\end{picture}\xrightarrow{\,\partial\,}
\Lambda^3\otimes\begin{picture}(42,12)(0,3)
\put(0,0){\line(1,0){30}}
\put(0,6){\line(1,0){42}}
\put(0,12){\line(1,0){42}}
\put(0,0){\line(0,1){12}}
\put(6,0){\line(0,1){12}}
\put(12,0){\line(0,1){12}}
\put(30,0){\line(0,1){12}}
\put(36,6){\line(0,1){6}}
\put(42,6){\line(0,1){6}}
\put(17,3){\makebox(0,0){.}}
\put(21,3){\makebox(0,0){.}}
\put(24,3){\makebox(0,0){.}}
\put(17,9){\makebox(0,0){.}}
\put(21,9){\makebox(0,0){.}}
\put(24,9){\makebox(0,0){.}}
\end{picture}\xrightarrow{\,\partial\,}\hspace{55pt}\\[6pt]
\cdots\xrightarrow{\,\partial\,}
\Lambda^\ell\otimes\begin{picture}(42,12)(0,3)
\put(0,0){\line(1,0){6}}
\put(0,6){\line(1,0){42}}
\put(0,12){\line(1,0){42}}
\put(0,0){\line(0,1){12}}
\put(6,0){\line(0,1){12}}
\put(12,6){\line(0,1){6}}
\put(30,6){\line(0,1){6}}
\put(36,6){\line(0,1){6}}
\put(42,6){\line(0,1){6}}
\put(17,9){\makebox(0,0){.}}
\put(21,9){\makebox(0,0){.}}
\put(24,9){\makebox(0,0){.}}
\end{picture}\xrightarrow{\,\partial\,}
\Lambda^{\ell+1}\otimes\begin{picture}(42,12)(0,3)
\put(0,6){\line(1,0){42}}
\put(0,12){\line(1,0){42}}
\put(0,6){\line(0,1){6}}
\put(6,6){\line(0,1){6}}
\put(12,6){\line(0,1){6}}
\put(30,6){\line(0,1){6}}
\put(36,6){\line(0,1){6}}
\put(42,6){\line(0,1){6}}
\put(17,9){\makebox(0,0){.}}
\put(21,9){\makebox(0,0){.}}
\put(24,9){\makebox(0,0){.}}
\end{picture}\to 0
\end{array}$$
is exact. Already the injectivity of the first homomorphism gives useful 
information regarding the higher Killing operator. Specifically it says that 
$$\begin{array}r
\begin{picture}(0,0)
\put(26,-6){\makebox(0,0){\tiny$\leftarrow\!\ell\!\rightarrow$}}
\put(49,-6){\makebox(0,0){\tiny$\leftarrow\!\ell\!\rightarrow$}}
\end{picture}
K_{apq\cdots rsbc\cdots de}
=K_{a(pq\cdots rs)(bc\cdots de)}\\[2pt]
K_{apq\cdots r(sbc\cdots de)}=0\\[2pt]
K_{[ap]q\cdots rsbc\cdots de}=0
\end{array}\hspace*{-5pt}\left.\rule{0pt}{19pt}\right\}
\enskip\Rightarrow\enskip K_{bpq\cdots rsbc\cdots de}=0.$$
In the flat case, if $X_{bc\cdots de}$ is a Killing tensor of valence $\ell$, 
it follows immediately from the Killing equation 
$\nabla_{(a}X_{bc\cdots de)}=0$, that 
$$K_{apq\cdots rsbd\cdots de}\equiv
\underbrace{\nabla_a\nabla_p\nabla_q\cdots\nabla_r\nabla_s}_{\ell+1}
X_{bc\cdots de}$$
satisfies exactly these symmetries and hence vanishes. In other words, the
Killing tensors of valence $\ell$ on ${\mathbb{R}}^n$ are polynomial of degree
at most~$\ell$. More generally, prolongation in the curved case implies that
the Killing tensors of valence $\ell$ near any point are determined by their
$\ell$-jet at that point.

\section{Lie algebra cohomology}
It remains to explain where (\ref{gen1}) and (\ref{gen2}) come from and the
answer is a special case of Kostant's generalised Bott-Borel-Weil
Theorem~\cite{k}, which we now explain. The special case we need involves only
the cohomology of an Abelian Lie algebra but for Kostant's results to apply it
is important that this Abelian Lie algebra be contained inside a semisimple Lie
algebra in a particular way. Specifically, let $${\mathfrak{g}}=
{\mathfrak{sl}}(n+1,{\mathbb{R}})=\{(n+1)\times(n+1)\mbox{ matrices }X
\mbox{ s.t.\ trace}(X)=0\}$$
and write ${\mathfrak{g}}=
{\mathfrak{g}}_{-1}\oplus{\mathfrak{g}}_0\oplus{\mathfrak{g}}_1$, comprising 
matrices of the form
$$\left\lgroup\begin{tabular}{c|ccc}
$0$&$0$&$\!\cdots\!$&$0$\\ \hline
$\ast$\\
\raisebox{2pt}[15pt]{$\vdots$}&&
\raisebox{-2pt}[0pt][0pt]{\makebox[0pt]{\Huge$0$}}\\
$\ast$\end{tabular}\right\rgroup\qquad
\left\lgroup\begin{tabular}{c|ccc}$\ast$&$0$&$\!\cdots\!$&$0$\\ \hline
$0$\\
\raisebox{2pt}[15pt]{$\vdots$}&&
\raisebox{-2pt}[0pt][0pt]{\makebox[0pt]{\Huge$\ast$}}\\
$0$\end{tabular}\right\rgroup\qquad
\left\lgroup\begin{tabular}{c|ccc}$0$&$\ast$&$\!\cdots\!$&$\ast$\\ \hline
$0$\\
\raisebox{2pt}[15pt]{$\vdots$}&&
\raisebox{-2pt}[0pt][0pt]{\makebox[0pt]{\Huge$0$}}\\
$0$\end{tabular}\right\rgroup,$$
respectively. Suppose ${\mathbb{V}}$ is an irreducible tensor representation
of~${\mathfrak{g}}$. It restricts to a representation of the Abelian subalgebra
${\mathfrak{g}}_{-1}$. Kostant's theorem computes the Lie algebra cohomology 
$H^p({\mathfrak{g}}_{-1},{\mathbb{V}})$. Explicitly, this means that the 
cohomology of the complex of ${\mathfrak{g}}_0$-modules
$$0\to{\mathbb{V}}
\xrightarrow{\,\partial\,}({\mathfrak{g}}_{-1})^*\otimes{\mathbb{V}}
\xrightarrow{\,\partial\,}\Lambda^2({\mathfrak{g}}_{-1})^*\otimes{\mathbb{V}}
\xrightarrow{\,\partial\,}\Lambda^3({\mathfrak{g}}_{-1})^*\otimes{\mathbb{V}}
\xrightarrow{\,\partial\,}\cdots$$
is computed as a ${\mathfrak{g}}_0$-module, where
$\partial:{\mathbb{V}}\to({\mathfrak{g}}_{-1})^*\otimes{\mathbb{V}}$ is defined
by the action of ${\mathfrak{g}}_{-1}$ on~${\mathbb{V}}$. To state the result,
we need a notation for the irreducible representations of 
${\mathfrak{sl}}(n+1,{\mathbb{R}})$ and for this we follow~\cite{be} writing, for 
example, 
$$\begin{picture}(90,5)
\put(0,0){\line(1,0){50}}
\put(70,0){\line(1,0){20}}
\put(0,0){\makebox(0,0){$\bullet$}}
\put(15,0){\makebox(0,0){$\bullet$}}
\put(30,0){\makebox(0,0){$\bullet$}}
\put(45,0){\makebox(0,0){$\bullet$}}
\put(60,0){\makebox(0,0){$\cdots$}}
\put(75,0){\makebox(0,0){$\bullet$}}
\put(90,0){\makebox(0,0){$\bullet$}}
\put(0,7){\makebox(0,0){\scriptsize$0$}}
\put(15,7){\makebox(0,0){\scriptsize$0$}}
\put(30,7){\makebox(0,0){\scriptsize$0$}}
\put(45,7){\makebox(0,0){\scriptsize$0$}}
\put(75,7){\makebox(0,0){\scriptsize$0$}}
\put(90,7){\makebox(0,0){\scriptsize$1$}}
\end{picture}\quad\mbox{and}\quad\begin{picture}(90,5)
\put(0,0){\line(1,0){50}}
\put(70,0){\line(1,0){20}}
\put(0,0){\makebox(0,0){$\bullet$}}
\put(15,0){\makebox(0,0){$\bullet$}}
\put(30,0){\makebox(0,0){$\bullet$}}
\put(45,0){\makebox(0,0){$\bullet$}}
\put(60,0){\makebox(0,0){$\cdots$}}
\put(75,0){\makebox(0,0){$\bullet$}}
\put(90,0){\makebox(0,0){$\bullet$}}
\put(0,7){\makebox(0,0){\scriptsize$1$}}
\put(15,7){\makebox(0,0){\scriptsize$0$}}
\put(30,7){\makebox(0,0){\scriptsize$0$}}
\put(45,7){\makebox(0,0){\scriptsize$0$}}
\put(75,7){\makebox(0,0){\scriptsize$0$}}
\put(90,7){\makebox(0,0){\scriptsize$0$}}
\end{picture}$$
for the defining representation ${\mathbb{R}}^{n+1}$ and its dual
$({\mathbb{R}}^{n+1})^*$, respectively. In particular, Kostant's theorem 
yields
$$H^0({\mathfrak{g}}_{-1},\;\begin{picture}(90,5)
\put(0,0){\line(1,0){50}}
\put(70,0){\line(1,0){20}}
\put(0,0){\makebox(0,0){$\bullet$}}
\put(15,0){\makebox(0,0){$\bullet$}}
\put(30,0){\makebox(0,0){$\bullet$}}
\put(45,0){\makebox(0,0){$\bullet$}}
\put(60,0){\makebox(0,0){$\cdots$}}
\put(75,0){\makebox(0,0){$\bullet$}}
\put(90,0){\makebox(0,0){$\bullet$}}
\put(0,7){\makebox(0,0){\scriptsize$0$}}
\put(15,7){\makebox(0,0){\scriptsize$\ell$}}
\put(30,7){\makebox(0,0){\scriptsize$0$}}
\put(45,7){\makebox(0,0){\scriptsize$0$}}
\put(75,7){\makebox(0,0){\scriptsize$0$}}
\put(90,7){\makebox(0,0){\scriptsize$0$}}
\end{picture}\;)=\;\begin{picture}(90,5)
\put(0,0){\line(1,0){50}}
\put(70,0){\line(1,0){20}}
\put(0,0){\makebox(0,0){$\times$}}
\put(15,0){\makebox(0,0){$\bullet$}}
\put(30,0){\makebox(0,0){$\bullet$}}
\put(45,0){\makebox(0,0){$\bullet$}}
\put(60,0){\makebox(0,0){$\cdots$}}
\put(75,0){\makebox(0,0){$\bullet$}}
\put(90,0){\makebox(0,0){$\bullet$}}
\put(0,7){\makebox(0,0){\scriptsize$0$}}
\put(15,7){\makebox(0,0){\scriptsize$\ell$}}
\put(30,7){\makebox(0,0){\scriptsize$0$}}
\put(45,7){\makebox(0,0){\scriptsize$0$}}
\put(75,7){\makebox(0,0){\scriptsize$0$}}
\put(90,7){\makebox(0,0){\scriptsize$0$}}
\end{picture}$$
where, again, we are following the \cite{be} to denote ${\mathfrak{g}}_0$ 
and its irreducible representations. More generally,
$$\begin{array}{rcl}
H^1({\mathfrak{g}}_{-1},\;\begin{picture}(90,5)
\put(0,0){\line(1,0){50}}
\put(70,0){\line(1,0){20}}
\put(0,0){\makebox(0,0){$\bullet$}}
\put(15,0){\makebox(0,0){$\bullet$}}
\put(30,0){\makebox(0,0){$\bullet$}}
\put(45,0){\makebox(0,0){$\bullet$}}
\put(60,0){\makebox(0,0){$\cdots$}}
\put(75,0){\makebox(0,0){$\bullet$}}
\put(90,0){\makebox(0,0){$\bullet$}}
\put(0,7){\makebox(0,0){\scriptsize$0$}}
\put(15,7){\makebox(0,0){\scriptsize$\ell$}}
\put(30,7){\makebox(0,0){\scriptsize$0$}}
\put(45,7){\makebox(0,0){\scriptsize$0$}}
\put(75,7){\makebox(0,0){\scriptsize$0$}}
\put(90,7){\makebox(0,0){\scriptsize$0$}}
\end{picture}\;)&=&\enskip\begin{picture}(90,5)
\put(0,0){\line(1,0){50}}
\put(70,0){\line(1,0){20}}
\put(0,0){\makebox(0,0){$\times$}}
\put(15,0){\makebox(0,0){$\bullet$}}
\put(30,0){\makebox(0,0){$\bullet$}}
\put(45,0){\makebox(0,0){$\bullet$}}
\put(60,0){\makebox(0,0){$\cdots$}}
\put(75,0){\makebox(0,0){$\bullet$}}
\put(90,0){\makebox(0,0){$\bullet$}}
\put(-3,7){\makebox(0,0){\scriptsize$-2$}}
\put(15,7){\makebox(0,0){\scriptsize$\ell\!+\!1$}}
\put(30,7){\makebox(0,0){\scriptsize$0$}}
\put(45,7){\makebox(0,0){\scriptsize$0$}}
\put(75,7){\makebox(0,0){\scriptsize$0$}}
\put(90,7){\makebox(0,0){\scriptsize$0$}}
\end{picture}\\[7pt]
H^2({\mathfrak{g}}_{-1},\;\begin{picture}(90,5)
\put(0,0){\line(1,0){50}}
\put(70,0){\line(1,0){20}}
\put(0,0){\makebox(0,0){$\bullet$}}
\put(15,0){\makebox(0,0){$\bullet$}}
\put(30,0){\makebox(0,0){$\bullet$}}
\put(45,0){\makebox(0,0){$\bullet$}}
\put(60,0){\makebox(0,0){$\cdots$}}
\put(75,0){\makebox(0,0){$\bullet$}}
\put(90,0){\makebox(0,0){$\bullet$}}
\put(0,7){\makebox(0,0){\scriptsize$0$}}
\put(15,7){\makebox(0,0){\scriptsize$\ell$}}
\put(30,7){\makebox(0,0){\scriptsize$0$}}
\put(45,7){\makebox(0,0){\scriptsize$0$}}
\put(75,7){\makebox(0,0){\scriptsize$0$}}
\put(90,7){\makebox(0,0){\scriptsize$0$}}
\end{picture}\;)&=&\hspace*{20pt}\begin{picture}(90,5)
\put(0,0){\line(1,0){50}}
\put(70,0){\line(1,0){20}}
\put(0,0){\makebox(0,0){$\times$}}
\put(15,0){\makebox(0,0){$\bullet$}}
\put(30,0){\makebox(0,0){$\bullet$}}
\put(45,0){\makebox(0,0){$\bullet$}}
\put(60,0){\makebox(0,0){$\cdots$}}
\put(75,0){\makebox(0,0){$\bullet$}}
\put(90,0){\makebox(0,0){$\bullet$}}
\put(-10,7){\makebox(0,0){\scriptsize$-\ell-3$}}
\put(15,7){\makebox(0,0){\scriptsize$0$}}
\put(30,7){\makebox(0,0){\scriptsize$\ell\!+\!1$}}
\put(45,7){\makebox(0,0){\scriptsize$0$}}
\put(75,7){\makebox(0,0){\scriptsize$0$}}
\put(90,7){\makebox(0,0){\scriptsize$0$}}
\end{picture}\\[7pt]
H^3({\mathfrak{g}}_{-1},\;\begin{picture}(90,5)
\put(0,0){\line(1,0){50}}
\put(70,0){\line(1,0){20}}
\put(0,0){\makebox(0,0){$\bullet$}}
\put(15,0){\makebox(0,0){$\bullet$}}
\put(30,0){\makebox(0,0){$\bullet$}}
\put(45,0){\makebox(0,0){$\bullet$}}
\put(60,0){\makebox(0,0){$\cdots$}}
\put(75,0){\makebox(0,0){$\bullet$}}
\put(90,0){\makebox(0,0){$\bullet$}}
\put(0,7){\makebox(0,0){\scriptsize$0$}}
\put(15,7){\makebox(0,0){\scriptsize$\ell$}}
\put(30,7){\makebox(0,0){\scriptsize$0$}}
\put(45,7){\makebox(0,0){\scriptsize$0$}}
\put(75,7){\makebox(0,0){\scriptsize$0$}}
\put(90,7){\makebox(0,0){\scriptsize$0$}}
\end{picture}\;)&=&\hspace*{20pt}\begin{picture}(90,5)
\put(0,0){\line(1,0){50}}
\put(70,0){\line(1,0){20}}
\put(0,0){\makebox(0,0){$\times$}}
\put(15,0){\makebox(0,0){$\bullet$}}
\put(30,0){\makebox(0,0){$\bullet$}}
\put(45,0){\makebox(0,0){$\bullet$}}
\put(60,0){\makebox(0,0){$\cdots$}}
\put(75,0){\makebox(0,0){$\bullet$}}
\put(90,0){\makebox(0,0){$\bullet$}}
\put(-10,7){\makebox(0,0){\scriptsize$-\ell-4$}}
\put(15,7){\makebox(0,0){\scriptsize$0$}}
\put(30,7){\makebox(0,0){\scriptsize$\ell$}}
\put(45,7){\makebox(0,0){\scriptsize$1$}}
\put(75,7){\makebox(0,0){\scriptsize$0$}}
\put(90,7){\makebox(0,0){\scriptsize$0$}}
\end{picture}\\[7pt]
\vdots\hspace*{40pt}&\vdots&\hspace*{40pt}\vdots\\[7pt]
H^{n-1}({\mathfrak{g}}_{-1},\;\begin{picture}(90,5)
\put(0,0){\line(1,0){50}}
\put(70,0){\line(1,0){20}}
\put(0,0){\makebox(0,0){$\bullet$}}
\put(15,0){\makebox(0,0){$\bullet$}}
\put(30,0){\makebox(0,0){$\bullet$}}
\put(45,0){\makebox(0,0){$\bullet$}}
\put(60,0){\makebox(0,0){$\cdots$}}
\put(75,0){\makebox(0,0){$\bullet$}}
\put(90,0){\makebox(0,0){$\bullet$}}
\put(0,7){\makebox(0,0){\scriptsize$0$}}
\put(15,7){\makebox(0,0){\scriptsize$\ell$}}
\put(30,7){\makebox(0,0){\scriptsize$0$}}
\put(45,7){\makebox(0,0){\scriptsize$0$}}
\put(75,7){\makebox(0,0){\scriptsize$0$}}
\put(90,7){\makebox(0,0){\scriptsize$0$}}
\end{picture}\;)&=&\hspace*{20pt}\begin{picture}(90,5)
\put(0,0){\line(1,0){50}}
\put(70,0){\line(1,0){20}}
\put(0,0){\makebox(0,0){$\times$}}
\put(15,0){\makebox(0,0){$\bullet$}}
\put(30,0){\makebox(0,0){$\bullet$}}
\put(45,0){\makebox(0,0){$\bullet$}}
\put(60,0){\makebox(0,0){$\cdots$}}
\put(75,0){\makebox(0,0){$\bullet$}}
\put(90,0){\makebox(0,0){$\bullet$}}
\put(-10,7){\makebox(0,0){\scriptsize$-\ell-n$}}
\put(15,7){\makebox(0,0){\scriptsize$0$}}
\put(30,7){\makebox(0,0){\scriptsize$\ell$}}
\put(45,7){\makebox(0,0){\scriptsize$0$}}
\put(75,7){\makebox(0,0){\scriptsize$0$}}
\put(90,7){\makebox(0,0){\scriptsize$1$}}
\end{picture}\\[7pt]
H^n({\mathfrak{g}}_{-1},\;\begin{picture}(90,5)
\put(0,0){\line(1,0){50}}
\put(70,0){\line(1,0){20}}
\put(0,0){\makebox(0,0){$\bullet$}}
\put(15,0){\makebox(0,0){$\bullet$}}
\put(30,0){\makebox(0,0){$\bullet$}}
\put(45,0){\makebox(0,0){$\bullet$}}
\put(60,0){\makebox(0,0){$\cdots$}}
\put(75,0){\makebox(0,0){$\bullet$}}
\put(90,0){\makebox(0,0){$\bullet$}}
\put(0,7){\makebox(0,0){\scriptsize$0$}}
\put(15,7){\makebox(0,0){\scriptsize$\ell$}}
\put(30,7){\makebox(0,0){\scriptsize$0$}}
\put(45,7){\makebox(0,0){\scriptsize$0$}}
\put(75,7){\makebox(0,0){\scriptsize$0$}}
\put(90,7){\makebox(0,0){\scriptsize$0$}}
\end{picture}\;)&=&\hspace*{35pt}\begin{picture}(90,5)
\put(0,0){\line(1,0){50}}
\put(70,0){\line(1,0){20}}
\put(0,0){\makebox(0,0){$\times$}}
\put(15,0){\makebox(0,0){$\bullet$}}
\put(30,0){\makebox(0,0){$\bullet$}}
\put(45,0){\makebox(0,0){$\bullet$}}
\put(60,0){\makebox(0,0){$\cdots$}}
\put(75,0){\makebox(0,0){$\bullet$}}
\put(90,0){\makebox(0,0){$\bullet$}}
\put(-16,7){\makebox(0,0){\scriptsize$-\ell-n-1$}}
\put(15,7){\makebox(0,0){\scriptsize$0$}}
\put(30,7){\makebox(0,0){\scriptsize$\ell$}}
\put(45,7){\makebox(0,0){\scriptsize$0$}}
\put(75,7){\makebox(0,0){\scriptsize$0$}}
\put(90,7){\makebox(0,0){\scriptsize$0$}}
\end{picture}\end{array}$$
where the right hand side follows the affine action of the Weyl group as 
explained in~\cite{be}. For our purposes, the crossed node can be dropped, 
viewing the results as irreducible tensor representations of 
${\mathfrak{sl}}(n,{\mathbb{R}})$. As tensor identities for 
${\mathfrak{sl}}(n,{\mathbb{R}})$, they are exactly what we need induce
(\ref{gen1}) and (\ref{gen2}) on a manifold.

\end{document}